\documentclass[11pt]{article} 

\usepackage[utf8]{inputenc} 
\usepackage{amssymb,amsmath}
\usepackage{amsfonts}

\usepackage{geometry} 
\usepackage{graphicx} 

\usepackage{booktabs} 
\usepackage{array} 
\usepackage{paralist} 
\usepackage{verbatim} 
\usepackage{subfig} 
\usepackage{ragged2e}
\usepackage{enumerate}
\usepackage{fancyhdr} 
\usepackage{sectsty}
\allsectionsfont{\sffamily\mdseries\upshape} 

\usepackage[nottoc,notlof,notlot]{tocbibind} 
\usepackage[titles,subfigure]{tocloft} 

\title{Finite Semihypergroups Built From Groups}
\author{Stan Onypchuk}

\begin{document}
\maketitle
\begin{abstract}
\begin{justify}
Necessary and sufficient conditions for finite semihypergroups to be built from groups of the same order are established.
\end{justify}
\end{abstract}

\noindent
\textbf{Introduction} The theory of hypergroups and semihypergroups was introduced by C. Dunkl [1], R. Jewett [2], and R. Spector [3] and is well developed now.\\
Many examples of finite commutative semihypergroups and hypergroups can be found in [4]. In [4] there is a precise \textsl{physical} definition of a finite semihypergroup - commutative semihypergroup because this example describes a finite collection of particles which interact by colliding.\\
Here is a more suitable example. Let us say we are observing a finite collection of events \{$e_{1},\ldots, e_{n}$\}, which may be in a cause and effect relationship. To establish that we will build a cube with frequencies (probabilities) $a_{i,j}(k)$ in which each event $e_{k}$ appears (be third) after any sequential pair of events $e_{i}$ and $e_{j}$ in our observed (infinitely long) sequence of events. So, we have a Markov chain of second order. Now assume that convolution - $e_{i} \ast e_{j} = \sum a_{i,j}(k)e_{k}$ - is associative and this semihypergroup (cube) can be developed from a group. In terms of this example the main result of present article shows that for such semihypergroups all observed events are not simple but rather are a combination of some simple events, evolution of which can be described as a group product for some appropriate group.\\
\\
Let $\mathbf{H}$ be a finite semihypergroup  with $n$ states $e_{1},e_{2}, \ldots , e_{n} $ and associative convolution operation defined by 
\begin{equation}
			e_{i} \ast e_{j} = \sum_{k} a_{i,j}(k)e_{k}  \text{  }    (i,j,k = 1,2, \ldots , n) 
\end{equation}			                      
where $a_{i,j}(k)\geqslant 0 \text{ and } \sum_{k = 1}^{n} a_{i,j}(k) = 1$  for each $i,j$.\\
Let us denote \\
columns \{$ a_{i,j}(1), \ldots , a_{i,j}(n)$\}  by $a_{i,j}$ \\
matrix with columns \{$a_{i,1}, \ldots , a_{i,n}$\}  by $A_{i} \text{ (matrix of left regular representation)} $\\
matrix with columns \{$a_{1,i}, \ldots , a_{n,i}$\}  by $B_{i} \text{ (matrix of right regular representation)}$ \\
and cube with matrices \{$A_{1}, \ldots , A_{n}$\} by $C.$ \\
\\
\textbf{Definition} A semihypergroup $\mathbf{H}$ will be said to be \emph{derived from a group} if any of the following two conditions is satisfied:\\
(A) There are only $n$ different rows in all $A_{i}$ and only $n$ different rows in all $B_{i}$ and in any $A_{i}$ rows are linearly independent and in any $B_{i}$ rows are linearly independent.\\
(B) There exists a group $\mathbf{G}$ of order $n$ with elements \{$g_{1}, \ldots , g_{n}$\} and measure $m$ on $\mathbf{G}$ that $e_{i} = m \cdot g_{i}$ where $g_{i} \cdot g_{j}$ is a group product and
\[
	e_{i} \ast e_{j} = (m \cdot g_{i}) \cdot (m \cdot g_{j})	
\]
\\
Now let us establish some properties which are shared among all semihypergroups with condition (A).\\
In matrix terms the convolution operation (1) transforms to $ a_{i,j} = A_{i}e_{j}$ (where $ e_{j} $ is a column-vector with 1 on $j$-th position and 0 otherwise). The matrix $\sum c_{i}A_{i} $ represents a measure $\sum c_{i}e_{i} $ so $\sum_{k=1}^{n} a_{i,j}(k)A_{k}$ represents $e_{i} \ast e_{j}$. Now it is obvious that the convolution operation (1) is associative if and only if
\begin{equation}
	A_{i}A_{j} = \sum_{k=1}^{n} a_{i,j}(k)A_{k} \\
\end{equation}
\textbf{Theorem 1.} If a semihypergroup satisfies condition (A) then there exist $n$ different non negative numbers among all $a_{i,j}(k)$ that $\sum a_{k} = 1$ and for each $a_{p}$ from \{$a_{1}, \ldots , a_{n}$\} and each $i, j = 1, \ldots , n$ there exist such $k$ and $l$ that 
\[
	a_{i,k}(j) = a_{p}  \text{ and } a_{l,i}(j) =a_{p}
\]
Let us denote as $r_{i,j} j$-th row in $A_{i}$. Then in $A_{1}A_{1}$ the first row is $\sum r_{1,1}(k)r_{1,k}$ and because of (2) it is  $\sum a_{1,1} r_{k,1}$. Now because $r_{1,k} (k = 1, \cdots ,n)$ are linearly independent we conclude that \{$r_{k,1}$\} contains the same set of linearly independent rows as \{$r_{1,k}$\} and the first column $a_{1,1}$ contains the same elements as the first row $r_{1,1}$.\\
The second row in $A_{1}A_{1}$ is $\sum r_{1,2}(k) r_{1,k} = \sum a_{1,1}(k) r_{k,2}$. Again we conclude that \{$r_{k,2}$\} contains the same set of linearly independent rows as \{$r_{1,k}$\} and the row $r_{1,2}$ contains the same elements as the column $a_{1,1}$.\\
Continuing this consideration for all $A_{i}A_{j}$, we conclude that each row in each $A_{i}$ and each column in $C$ contains the same  elements. Now, if we repeat this process with matrices $B_{i}$, we will see that each row in each matrix $B_{i}$ contains the same elements as each column in $C$.\\ 
Now let us establish some properties which are shared among all semihypergroups derived from a group.\\
\textbf{Corollary 1.} Diagonal elements in each matrix $A_{i}$ and $B_{i}$ are equal.\\
The $i$-th row in  $A_{1}A_{1}$ is
\[
	\sum r_{1,i}(k)r_{1,k} =\sum  a_{1,1}r_{k,i}
\] 
when here in the left part $k = i$ and in the right part $k = 1$ we have $r_{1,i}(i) = a_{1,1}(1) = r_{1,1}(1)$.\\ 
\textbf{Corollary 2.} Every matrix $A_{i}$ is a linear combination of $G_{i}$, where \{$G_{i}$\} - matrices of a left regular representation of some group $\mathbf{G}$ of order $n$.\\
To prove this we will use the following alternative definition of a group:\\
\\
\hspace*{0.3in}\textbf{a)} On non-empty set $G$ of elements \{$g_{i}$\} with a binary operation - $g_{i} \cdot g_{j} = g_{k}$;\\
\hspace*{0.3in}\textbf{b)} This operation is associative;\\
\hspace*{0.3in}\textbf{c)} For any elements $g_{i}$ and $g_{j}$ there exist at least one such $g_{k}$ and at least one such $g_{l}$ that
\centerline{	$g_{k} \cdot g_{i} = g_{i} \cdot g_{l} = g_{j}.$}
\\
By condition (B) there are no more than $n$ different elements in $C$. Let us assume that all $n$ elements are different and one of them is $a$.\\
Then for any $i,k$ there exists such $j$ that
\begin{equation}
	a_{i,j}(k) = a
\end{equation}
and for any $j,k$ there exists such $i$ that
\begin{equation}
	a_{i,j}(k) = a
\end{equation}
When $a$ = 1 and all other elements in the cube $C$ are equal to 0, conditions (3) and (4) describe binary operation defined as
\begin{equation}
	g_{i} \cdot g_{j} = a_{i,j}(k) = g_{k}
\end{equation}
Let us show that the operation (5) is associative.
\begin{equation}
	(e_{i} \ast e_{j}) \ast e_{m} = \sum_{k=1}^{n} a_{i,j}(k)(e_{k} \ast e_{m}) = \sum_{k=1}^{n}a_{i,j}(k)a_{k,m}
\end{equation}
On the other hand,
\begin{equation}
	e_{i} \ast (e_{j} \ast e_{m}) = \sum_{p=1}^{n} a_{j,m}(p)(e_{i} \ast e_{p}) = \sum_{p=1}^{n}a_{j,m}(p)a_{i,p}
\end{equation}
Now from (6) and (7) 
\[
	\sum a_{i,j}(k) a_{k,m} = \sum a_{j,m}(p)a_{i,p}
\]
and because $a_{i,j}(k) = a_{j,m}(p) = 1$ that means that
\[
	(g_{i} \cdot g_{j}) \cdot g_{m} = g_{i} \cdot (g_{j} \cdot g_{m}).
\]
Therefore, a cube $C$ with only one element $a_{p}$ equal to 1 and all other elements equals to 0, describes a group.\\
Let us show that different $a_{i}, a_{j}$ describe the same group. Assume that $a_{i}$ and $a_{j}$ correspondent to two different groups with its left representations $G_{1}$ and $G_{2}$.\\
Consider the following two fusions of their matrices \\
\begin{equation}
	G_{1,i_{1}}, \ldots , G_{1,i_{k}}, \ldots , G_{1,i_{n}} 
\end{equation}
\begin{equation}
	G_{1,j_{1}}, \ldots , G_{1,j_{k}}, \ldots , G_{1,j_{n}} 
\end{equation}
\begin{equation}
	G_{2,j_{1}}, \ldots , G_{2,j_{k}}, \ldots , G_{2,j_{n}}
\end{equation}

Here in (8) and (9) we have matrices from representation of the same group in different order when in pair (9) and (10) there are matrices from representation of different groups and $G_{1,j_{k}} \ne G_{2,j_{k}}$. Then (let us assume there are only two numbers $a_{1}$ and $a_{2}$ not equal to zero) if matrices in (8) and (9) are in such order that every matrix
\[
a_{1}G_{1,i_{1}} + a_{2}G_{1,j_{1}}, \ldots , a_{1}G_{1,i_{k}} + a_{2}G_{1,j_{k}}, \ldots , a_{1}G_{1,i_{n}} + a_{2}G_{1,j_{n}}
\]
has the same rows then in combination (8) and (10) with the same coefficients $a_{1}$ and $a_{2}$ there will be a matrix with the same rows as for (8) and (9) (when $G_{2,j_{j}}$ is an identity matrix, for example) and there will be a matrix with different set of rows when $G_{2,j_{i}}$ is not equal to any of matrices from (9). This contradicts condition (A).\\
\textbf{Theorem 2}. If associative semihypergroup \textbf{H} satisfies condition (A), then it satisfies condition (B).\\
It has been shown above that every matrix $A_{i}$ of left regular representation of \textbf{H} is a combination of $G_{i}$ - matrices of left representation of some group \textbf{G}. By (A) any $A_{i}$ should have the same set of rows that the matrix $A_{1}$ has. Let us assume that
\[
	A_{1} = a_{1}G_{1} + a_{2}G_{2} + \ldots + a_{n}G_{n}
\]
In the product $G_{i}A_{1}$, matrix $G_{i}$ acts as a permutation of rows in matrix $A_{1}$. So, the only possibility to satisfy the following two conditions: $A_{i}$ as a combination of $G_{i}$ and $A_{i}$ contains the same set of rows as $A_{1}$ is
\[
	A_{i} = G_{i}A_{1} = a_{1}G_{i}G_{1} + a_{2}G_{i}G_{2} + \ldots + a_{n}G_{i}G_{n}
\]
The first means that $A_{i}$ contains the same rows as $A_{1}$ and the second that $A_{i}$ is a combination of $G_{i}$.\\
The rest should be clear from the following consideration\\
\textbf{Theorem 3}. If an associative semihypergroup \textbf{H} satisfies condition (B) then it satisfies condition (A).\\
Let \textbf{G} be a group of order $n$ and $g_{i} \cdot g_{j}$ denote product operation and $m$ a probability measure on \textbf{G}. Let us build a semihypergroup \textbf{H} with elements \{$e_{1}, \ldots , e_{n}$\} where $e_{i} = m \cdot g_{i}$ and convolution operation as 
\[
	e_{i} \ast e_{j} = (m \cdot g_{i}) \cdot (m \cdot g_{j})
\] 
Because elements of \textbf{H} are $m \cdot g_{i}$, the convolution should be described in terms of $(m \cdot g_{i})$.
\[
	e_{i} \ast e_{j} = m \cdot g_{i} \cdot m \cdot g_{j} =  m \cdot (g_{i} \cdot m \cdot g_{j}) = m\cdot (\sum m_{i,j}(k)g_{k})
\]
Let us consider how the $i$-th matrix $M_{i}$ of a regular representation of $\mathbf{H}$ is built.\\
Each $m_{i,j}$ the $j$-th column in $M_{i}$ is the ordered sequence of coefficients in $g_{i}\cdot m \cdot g_{j}$ - \{$m_{i,j}(k)$\}.\\
Because 
\[ 
	m \cdot g_{j} = m_{1}g_{1} \cdot g_{j} + m_{2}g_{2} \cdot g_{j} + \ldots + m_{n}g_{n} \cdot g_{j}
\]
and $g_{i} \cdot g_{j}$ is the $j$-th column in $G_{i}$ (more precisely $j$-th column in $G_{i}$ is the sequence \{$g_{i,j}(1), \ldots , g_{i,j}(n)$\}  of coefficients in $g_{i} \cdot  g_{j} = g_{i,j}(1)g_{1} + g_{i,j}(2)g_{2} + \ldots + g_{i,j}(n)g_{n}$) we have to conclude that $m \cdot g_{j}$ represents the column $\sum m_{k} g_{k,j}$ and $m\cdot g_{j}  (j = 1, \ldots , n)$ should be represented by matrix 
\[
\sum m_{k} G_{k} = M = 
\begin{pmatrix}
		m_{1,1}(1) & m_{1,2}(1) & \ldots & m_{1,n}(1) \\
		m_{1,1}(2) & m_{1,2}(2) & \ldots & m_{1,n}(2) \\
		\ldots & \ldots & \ldots \\
		\ldots & \ldots & \ldots \\
		\ldots & \ldots & \ldots \\
		m_{1,1}(n) & m_{1,2}(n) & \ldots & m_{1,n}(n)
	\end{pmatrix}
\]
where $m_{1,1}(1) = m_{1}, m_{1,1}(2) = m_{2}, \ldots ,m_{1,1}(n) = m_{n}$.
Now to add scalar multiplier $g_{i}$ to $M$, we need to implement the multiplication by component $g_{i}$ on $M$. It means that each $k$-th row in $M$ ($m_{1,1}(k), m_{1,2}(k), \ldots ,m_{1,n}(k)$) that is \{$m_{1,1}(k)g_{i}\cdot g_{k}, m_{1,2}(k)g_{i}\cdot g_{k}, \ldots ,m_{1,n}(k)g_{i}\cdot g_{k}$\} will be moved to $p$-th row where $g_{p} = g_{i} \cdot g_{k}$.\\
It means that $M_{i} = G_{i} M$. Because in $G_{i}M$,  $ G_{i}$ acts as permutation of rows in $M$, any $M_{i}$ has the same set of rows as $M$.\\  
The same is true with respect to the right regular representation of \textbf{H}. Next, if measure $m$ is not a uniform measure concentrated on some subgroup of \textbf{G}, all rows in $M_{i}$ will be different and linearly independent.\\
Semihypergroup \textbf{H} is associative because underlying group \textbf{G} is associative and $e_{i} \ast e_{j} = (m \cdot g_{i}) \cdot (m \cdot g_{j}) = m \cdot (g_{i} \cdot m \cdot g_{j})$. \\ 
\\
\textbf{References}\\
\text{[1]} C. F. Dunkl, The measure algebra of a locally compact hypergroup, Trans. Amer. Math. Soc. 179(1973) 331-348.\\
\text{[2]} R. I. Jewett, Spaces with an abstract convolution of measures, Advances in Math. 18 (1975) 1-110.\\
\text{[3]} R. Spector, Mesures invariantes sur les hypergroupes,Trans Amer Math Soc 239, (1978) 147 -n 165.\\
\text{[4]} N. J. Wildberger, Finite commutative hypergroups and applications from group theory to conformal field theory, Contemporary Mathematics, 188 (1995) 413-431. 

\end{document}